\documentclass[11pt,oneside]{amsproc}

\usepackage[cp1251]{inputenc}
\usepackage[russian]{babel}
\usepackage{amssymb,epic,eepic}
\numberwithin{equation}{section}
\newtheorem{lem}{Лемма}[section]
\newtheorem{tm}{Теорема}[section]

\newtheorem{rem}{Замечание}[section]
\newtheorem{dif}{Определение}[section]

\title[]{Особые точки самоподобной функции нулевого спектрального порядка. Самоподобная струна Стилтьеса}
\author{И.~А.~Шейпак\footnote{%
Работа выполнена при поддержке грантов РФФИ \No\,07-01-00283,  поддержки ведущих научных школ \No\,НШ-2513.2008.1 и фондом INTAS, грант №
05-1000008-7883}
\address{Московский государственный университет
им.~М.~В.~Ломоносова, механико-математический факультет}
\email{iasheip@mech.math.msu.su}}
\newcommand{\Wo}{{\raisebox{0.2ex}{\(\stackrel{\circ}{W}\)}}{}}

\newcommand{\sign}{\operatorname{sign}}

\begin{document}
\noindent УДК~517.518,517.984
\begin{abstract}
В статье вводиться понятие самоподобной функции нулевого спектрального порядка и  изучаются её свойства. Эта функция
имеет не более чем счётное число точек разрыва, причём все точки разрыва являются точками разрыва 1-го рода, за
исключением быть может одной точки, являющейся \emph{особой}. Получена формула, позволяющая по параметрам самоподобия
функции, вычислить её координаты. Исследуется поведение самоподобной функции в окрестности особой точки.

Неубывающая функция $f$ нулевого спектрального порядка, принадлежащая пространству $L_2[0,1]$, порождает самоподобную
струну Стилтьеса, т.е. спектральную задачу вида
\begin{gather*}
    -y''-\lambda\rho y=0,\\ y(0)=y(1)=0,
\end{gather*}
где $\rho$ есть функция из пространства $\Wo_2^{-1}[0,1]$ и $f'=\rho$. Незнакоопределённая функция $f$ приводит к
понятию самоподобной индефинитной струны Стилтьеса.
\end{abstract}
\maketitle

\section{Введение}

Струна называется стилтьесовской, если она представляет собой невесомую нить, которая несет лишь сосредоточенные массы
$m_1$, $m_2,$ $\ldots$ в точках $0=x_1<x_2<\ldots$, сгущающихся к правому концу, $\lim x_n=L(\leqslant \infty)$. Задача
о колебаниях стилтьесовской струны единичной длины c закрепленными концами приводит к спектральной задаче
\begin{gather}\label{eq:vved1}
-y''-\lambda\rho y=0,\\
y(0)=y(1)=0,\label{eq:vved2}
\end{gather}
где $\rho=\sum_{k=1}^\infty m_k \delta(x-x_k)$. Обобщённая первообразная $f$ ($f'=\rho$) такого веса представляет собой кусочно-постоянную
функцию, точки разрыва которой стремятся к правому концу отрезка, а $f(1-0)=\sum_{k=1}^n m_k (\leqslant \infty)$. C этой точки зрения, можно
считать, что точка $x=1$ является особой точкой функции $f$.

Более подробно о стилтьесовской струне см. в~\cite{KacKr}.

Укажем на более широкий  класс спектральных задач, связанных со стилтьесовской струной. В работе~\cite{SV} изучалась
задача~\eqref{eq:vved1}--\eqref{eq:vved2} в предположении, что вес $\rho$ является самоподобной мерой. Но случай
дискретных, т.\,е. как раз случай стилтьесовской струны не рассматривается. Для таких классов весов были получены
асимптотические формулы распределения собственных значений, имеющие степенной порядок.

Самоподобные функции положительного спектрального порядка были введены в работах~\cite{VlSh1}--\cite{VlSh2}.
Оказывается, что спектральный порядок функции $f$ непосредственно связан с порядком асимптотики собственных значений
задачи~\eqref{eq:vved1}--\eqref{eq:vved2}, где
$\rho=f'$. В этих работах изучалась задача~\eqref{eq:vved1}--\eqref{eq:vved2} с более общим классом весовых функций. А
именно, в качестве веса $\rho$ рассматривались самоподобные обобщённые функции из пространства
\(\Wo_2^{-1}[0,1]\).  Первообразные таких весов являются самоподобными
функциями из пространства $L_2[0,1]$. При этом условия, накладываемые на параметры самоподобия исключают из
рассмотрения кусочно-постоянные функции, т.\,е. и в этих работах стилтьесовские струны не рассматривались. Заметим, что
собственные значения задачи~\eqref{eq:vved1}--\eqref{eq:vved2} с весом, являющимся обобщённой производной функции
положительного спектрального порядка также подчиняются степенному закону распределения.

С точки зрения задачи~\eqref{eq:vved1}--\eqref{eq:vved2} случай, когда вес является обобщённой производной самоподобной
функции нулевого спектрального порядка занимает особое место. А именно, собственные значения имеют экспоненциальное
распределение (см.\cite{VlSh3}). Поэтому изучение таких функций представляет интерес с точки зрения спектральной теории
операторов.

Отметим, что самоподобные функции (в том числе и кусочно-постоянные) имеют важное значение не только для спектральных задач. Например, одним из
популярных приложений самоподобных непрерывных функций является теория сжатия компьютерных изображений~(см. \cite{Bar}). Различные задачи,
связанные с изучением таких свойств самоподобных функций как ограниченность вариации, абсолютная непрерывность и других свойств, естественным
образом приводящих к понятию  вероятностных самоподобных мер, возникают в теории фрактальных кривых, изучение которых началось с уже ставших
классическими работы \cite{deRam1}--\cite{deRam3}. Среди более поздних работ, посвящённых различным свойствам фрактальных кривых,  отметим
работы~\cite{Prot1} и~\cite{Prot2}. Конструкция аффинно-самоподобных функций, принадлежащих пространствам $L_p[0,1]$, приведена в~\cite{Sheip}.
Кроме того, самоподобные меры находят применение и в теории случайных процессов (см. например, \cite{ANaz}).

Целью статьи является определение функции нулевого спектрального порядка, исследование её структуры, классификация её
особых точек в зависимости от параметров самоподобия и исследование поведения самоподобной функции в окрестности этой
точки.

Также вводится понятие самоподобной струны Стилтьеса, имеющей тесную связь с функциями нулевого спектрального порядка. В отличие от классической
струны функции нулевого спектрального порядка порождают вообще говоря индефинитную струну (ср.~\cite{KacKr}, замечание~\ref{rem:indef_string},
теорема~\ref{tm:def_string})

Напомним основные понятия, связанные с самоподобными функциями.

\section{Самоподобные функции в пространстве $L_p[0,1]$ и их спектральный порядок}\label{pt:666}
\subsection{Операторы подобия в пространстве $L_p[0,1]$}
Пусть фиксировано натуральное число $n>1$, и пусть вещественные числа
$a_k>0$, где \(k=1,\ldots,n\), таковы, что
$$
\sum\limits_{k=1}^n a_k=1,
$$
а вещественные числа $d_k$ и $\beta_k$ пока произвольны. Дополнительно определим числа $\alpha_1=0$,
$\alpha_k=\sum_{l=1}^{k-1}a_l$, $k=2,\ldots,n+1$. Таким образом, числа
$0=\alpha_1<\alpha_2<\ldots<\alpha_n<1=\alpha_{n+1}$ образуют разбиение отрезка $[0,1]$.

Определим непрерывные операторы $G_k$ и $\tilde G_k$, сопоставляющие функции $f\in L_p[0,1]$ функцию из $L_p[\alpha_k,\alpha_{k+1}]$,
$k=1,2,\ldots,n$ по правилу:
\begin{gather*}
G_k(f)(x)=f(t), \quad \text{где} \quad x=a_kt+\alpha_k, \quad x\in [\alpha_k,\alpha_{k+1}],\quad t\in[0,1]; \\
\tilde G_k(f)(x)=f(t), \quad \text{где} \quad x=-a_kt+\alpha_{k+1}, \quad x\in [\alpha_k,\alpha_{k+1}],\quad t\in[0,1].
\end{gather*}

Окончательно построим нелинейный непрерывный оператор $G:L_p[0,1]\to L_p[0,1]$ вида
\begin{equation}\label{eq:auxto}
G(f)=\sum\limits_{k=1}^n\left\{\beta_k\cdot\chi_{(\alpha_k,\alpha_{k+1})} +d_k\cdot \hat G_k(f)\right\},
\end{equation}
где $\chi_{(\alpha_k,\alpha_{k+1})}$ --- характеристическая функция интервала $(\alpha_k,\alpha_{k+1})$, а
$\hat G_k=G_k$ или $\hat G_k=\tilde G_k$ в зависимости от $k$.
Про отображения $\tilde G_k$ будем говорить, что они меняют ориентацию отрезка $[\alpha_k,\alpha_{k+1}]$.

Операторы $G$ вида~\eqref{eq:auxto} будут называться \emph{операторами подобия}.

В работе~\cite{Sheip} рассмотрены более общие операторы подобия.
\begin{equation}\label{eq:oper_pod_ck}
G(f)=\sum\limits_{k=1}^n\left\{\beta_k\cdot\chi_{(\alpha_k,\alpha_{k+1})}+c_k\cdot x +d_k\cdot G_k(f)\right\},
\end{equation}
но не рассматривались операторы, меняющие ориентацию отрезков $[\alpha_k,\alpha_{k+1}]$.

Там же доказана следующая
\begin{lem}\label{lem2:1}
Оператор подобия $G$ является сжимающим в $L_p[0,1]$ в том и только том случае, когда справедливы неравенства
\begin{gather}\label{eq:szim1}
\sum\limits_{k=1}^n a_k\,|d_k|^p<1 \quad (1\leqslant p<+\infty),\\\label{eq:szim2}
\max_{1\leqslant k\leqslant n}|d_k|<1 \quad (p=+\infty).
\end{gather}
\end{lem}

Из леммы~\ref{lem2:1} и принципа сжимающих отображений немедленно вытекает справедливость
следующего утверждения.
\begin{tm}\label{sek2:1}
Если справедливо неравенство~\eqref{eq:szim1}(\eqref{eq:szim2}), то существует и единственна функция $f\in L_p[0,1]$,
удовлетворяющая уравнению $G(f)=f$.
\end{tm}

В дальнейшем всегда будет предполагаться, что либо неравенство~\eqref{eq:szim1} либо неравенство~\eqref{eq:szim2}
выполнено.

Если функция $f\in L_p[0,1]$ удовлетворяет уравнению $G(f)=f$, где $G$ --- некоторый оператор подобия, то такая функция
будет называться \emph{самоподобной}. При этом величины $n$, $a_k$, $d_k$ и $\beta_k$, где $k=1,2,\ldots,n$,
определяющие соответствующий оператор подобия \(G\), будут называться \emph{параметрами самоподобия} функции $f$.

\subsection{Спектральный порядок}
Заметим, что если для всех $k=1,2,\ldots,n$ выполнено
$\beta_k=0$, то оператор подобия $G$ имеет только тривиальную неподвижную точку $f\equiv 0$,
поэтому в дальнейшем будем предполагать, что выполнено условие
\begin{enumerate}
\item[($B$)] среди чисел \(\beta_k\), где \(k=1,\ldots,n\), по меньшей мере одно отлично от нуля.
\end{enumerate}

Среди самоподобных функций, удовлетворяющих условию $B$, выделим  следующие классы, для которых
параметры самоподобия соответственно обладают свойствами:
\begin{enumerate}
\item[($D_0$)] \(d_k=0\) для всех \(k=1,\ldots,n\);
\item[($D_1$)] среди чисел \(d_k\), где \(k=1,\ldots,n\), ровно одно отлично от нуля;
\item[($D_2$)] среди чисел \(d_k\), где \(k=1,\ldots,n\), не менее двух отличны от нуля;
\end{enumerate}

Самоподобные функции, принадлежащие пространству $L_2[0,1]$ и  удовлетворяющие условию $D_2$, были введены в
работах~\cite{VlSh1},\cite{VlSh2} и называются
\emph{самоподобными функциями положительного спектрального порядка}. В этой же работе показано, что для параметров
самоподобия, отвечающих такой функции, уравнение
\begin{equation}\label{eq:spord}
\sum\limits_{k=1}^n\left(a_k\,|d_k|\right)^{D}=1
\end{equation}
имеет решение \(0<D<1\), а  для собственных значений $\lambda_n$ задачи~\eqref{eq:vved1}--\eqref{eq:vved2} справедлива асимптотическая формула
$$
|\lambda_n|\asymp n^{1/D}.
$$

Нетрудно заметить, что самоподобные функции, удовлетворяющие условию $D_0$, суть простые функции, принимающие конечное
число значений, а именно
$$
f=\sum_{k=1}^n \beta_{k}\cdot\chi_{(\alpha_k,\alpha_{k+1})}.
$$
Обобщенная производная такой функции имеет вид
\begin{equation}\label{eq:konech}
 f'(x)=\sum_{k=2}^n
(\beta_{k}-\beta_{k-1})\cdot\delta(x-\alpha_k).
\end{equation}

Подстановка~\eqref{eq:konech}  в качестве веса $\rho$ в~\eqref{eq:vved1}--\eqref{eq:vved2} сводит задачу к
конечномерной. Непосредственное вычисление позволяет убедиться, что собственные значения в этом случае являются корнями
многочлена
\begin{multline*}
1-\lambda\sum_{k=2}^n m_k \alpha_k(1-\alpha_k)+\lambda^2\sum_{2\leqslant i_1<i_2\leqslant n}
m_{i_1}m_{i_2}\alpha_{i_1}(\alpha_{i_2}-\alpha_{i_1})(1-\alpha_{i_2})+\ldots+\\
+(-1)^k\lambda^k \sum_{2\leqslant i_1<i_2<\ldots<i_k\leqslant n}m_{i_1}m_{i_2}\ldots m_{i_k}
\alpha_{i_1}(\alpha_{i_2}-\alpha_{i_1})\ldots(\alpha_{i_k}-\alpha_{i_{k-1}})(1-\alpha_{i_2})+\ldots\\
\ldots+ (-1)^{n-1} \lambda^{n-1} m_{2}m_{3}\ldots
m_{n}\alpha_{2}(\alpha_{3}-\alpha_{2})\ldots(\alpha_{n}-\alpha_{n-1})(1-\alpha_{n}),
\end{multline*}
где $m_k=\beta_k-\beta_{k-1}$, $k=2,3,\ldots,n$.

Для самоподобных функций класса $D_1$ уравнение~\eqref{eq:spord} имеет  только тривиальное решение
$D=0$. В связи с этим дадим следующее
\begin{dif} Самоподобные функции класса $D_1$ будем называть самоподобными функциями нулевого спектрального
порядка.
\end{dif}

Изучим более подробно функции класса $D_1$ (не обязательно принадлежащие пространству $L_2[0,1]$). Обозначим через
$\hat k$ тот единственный индекс, для которого
$d_{\hat k}\ne 0$, $1\leqslant \hat k\leqslant n$.

Справедливо следующее утверждение.
\begin{tm}\label{tm:schet}
Самоподобная функция нулевого спектрального порядка является кусочно-постоянной и принимает не более чем счётное число
значений. Все точки разрыва являются точками разрыва 1-го рода, кроме быть может одной точки.
\end{tm}
\begin{proof}
Рассмотрим сначала случай, когда в формуле~\eqref{eq:auxto} все операторы $\hat G_k$ равны $G_k$. Случай, когда при некоторых
$k\in[1,n]$ оператор $\hat G_k=\tilde G_k$, рассматривается аналогично. Функция $f$, являющаяся неподвижной точкой оператора $G$, в этом случае
имеет некоторые отличия, о которых мы расскажем в \S 3.

Для каждого натурального $m$ построим разбиение отрезка $[0,1]$ на $n^m$ отрезков и каждому такому отрезку разбиения сопоставим
последовательность чисел $k_1,k_2,\ldots, k_m$ (номер отрезка), где
$k_i=1,2,\ldots,n$, $i=1,2,\ldots,m$.

При $m=1$ такое разбиение образуют отрезки $I_{k_1}\colon=[\alpha_{k_1},\alpha_{k_1+1}]$,
$k_1=1,2,\ldots,n$. Такие отрезки будем называть \emph{отрезками разбиения первого порядка} или
просто \emph{отрезками первого порядка}.

При $m=2$ разбиение образуют отрезки $I_{k_1,k_2}\colon=\{y=a_{k_2}x+\alpha_{k_2},\quad x\in
I_{k_1}\}$, $k_i=1,2,\ldots,n$, $i=1,2$. Такие отрезки будем называть \emph{отрезками разбиения
второго порядка} или просто \emph{отрезками второго порядка}.

Далее построение осуществляется индуктивно. Если построено разбиение на отрезки
$I_{k_1,k_2,\ldots,k_{m-1}}$, $k_i=1,2,\ldots, n$, $i=1,2,\ldots, m-1$, то очередное разбиение
определяется следующим образом $I_{k_1,k_2,\ldots,k_{m}}\colon=\{y=a_{k_m}x+\alpha_{k_m},\quad x\in
I_{k_1,k_2,\ldots,k_{m-1}}\}$, $k_i=1,2,\ldots,n$, $i=1,2,\ldots,m$. Такие отрезки будем называть
\emph{отрезками разбиения $m$-го порядка} или просто \emph{отрезками $m$-го порядка}. Нетрудно
заметить, что отрезок $I_{k_1,k_2,\ldots,k_{m}}$ $m$-го порядка имеет длину $a_{k_1}a_{k_2}\ldots
a_{k_m}$.

Пусть $f_0\equiv 0$. Построим последовательность $f_j$, $j=1,2,\ldots$ по правилу
$f_{j}=G(f_{j-1})$. Очевидно, что $f_j\to f$ в $L_p[0,1]$. При этом в силу
соотношение~\eqref{eq:auxto} выполнено

\begin{align}
\label{eq:f1}
&f_1(x)=\beta_{k_1}\text{ при } x\in I_{k_1},\quad k_1=1,2,\ldots,n,\\
\label{eq:f2}
&f_2(x)=\left\{\begin{aligned}&\beta_{k_1}, \quad x\in I_{k_1}, \quad k_1\ne \hat k,\\
&d_{\hat k}\beta_{k_2}+\beta_{\hat k},\quad x\in I_{\hat k,k_2},\quad k_2=1,2,\ldots, n,
\end{aligned}\right.\\
&f_3(x)=\left\{\begin{aligned}&\beta_{k_1}, \quad x\in I_{k_1}, \quad k_1\ne \hat k,\\
&d_{\hat k}\beta_{k_2}+\beta_{\hat k},\quad x\in I_{\hat k,k_2},\quad k_2\ne\hat k,\\
&d_{\hat k}(d_{\hat k}\beta_{k_3}+\beta_{\hat k})+\beta_{\hat k},\quad, x\in I_{\hat k,\hat k,k_3},
\quad k_3=1,2,\ldots,n\\
\end{aligned}\right.\\
\notag
&\ldots\\
\label{eq:fm}
&f_m(x)=\left\{\begin{aligned}& \beta_{k_1}, \quad x\in I_{k_1}, \quad k_1\ne \hat k,\\
&d_{\hat k}\beta_{k_2}+\beta_{\hat k},\quad x\in I_{\hat k,k_2},\quad k_2\ne \hat k,\\
&\ldots\\
&d_{\hat k}^{m-2}\beta_{k_{m-1}}+\beta_{\hat k}(1+d_{\hat k}+\ldots+d_{\hat k}^{m-3}),\quad
x\in I_{\tiny\underbrace{\hat k,\hat k,\ldots,\hat k}_{m-2},k_{m-1}},\quad k_{m-1}\ne \hat k\\
&d_{\hat k}^{m-1}\beta_{k_m}+\beta_{\hat k}(1+d_{\hat k}+\ldots+d_{\hat k}^{m-2}),\quad x\in I_{\tiny\underbrace{\hat
k,\hat k,\ldots,\hat k}_{m-1},k_m},
\quad k_m=1,2,\ldots, n.
\end{aligned}\right.
\end{align}

Несложно заметить, что функции $f_m$ обладают тем свойством, что вне отрезка $I_{\underbrace{\hat k,\hat k,\ldots,\hat
k}_{m}}$
$$
G(f_m)=f_m,
$$
где последнее равенство выполнено в смысле пространства $L_p[0,1]$. Следовательно, вследствие самоподобия функции $f$
вне отрезка $I_{\underbrace{\hat k,\hat k,\ldots,\hat k}_{m}}$ выполнены равенства
$$
f_m=f.
$$
Таким образом, функции $f_m$ равны функции $f$ вне одного отрезка $m$-го порядка.

Отрезки $I_{\underbrace{\hat k,\hat k,\ldots,\hat k}_{m}}$ образуют систему вложенных отрезков, стягивающихся к  точке
$\hat x$, являющейся в некотором роде особой для самоподобной функции $f$. Левый ($x_l$) и правый ($x_r$) концы отрезка
$I_{\underbrace{\hat k,\hat k,\ldots,\hat k}_{m}}$ соответственно имеют координаты
$$
x_{l}=\alpha_{\hat k}+a_{\hat k}\alpha_{\hat k}+a^2_{\hat k}\alpha_{\hat k}+\ldots+a^{m-1}_{\hat k}\alpha_{\hat
k},\quad x_{r}=\alpha_{\hat k}+a_{\hat k}\alpha_{\hat k}+a^2_{\hat k}\alpha_{\hat k}+\ldots+a^{m-2}_{\hat
k}\alpha_{\hat k}+a^{m-1}_{\hat k}\alpha_{\hat k+1}.
$$

Координату точки $\hat x$ можно вычислить по формуле
\begin{equation}\label{eq:osob_tochka}
\hat x=\dfrac{\alpha_{\hat k}}{1-a_{\hat k}}.
\end{equation}
В частности, если $\hat k=n$, т.\,е. ровно последнее $d_k$ отлично от нуля, то из~\eqref{eq:osob_tochka} и определений
чисел $\{a_k\}$, $\{\alpha_k\}$ вытекает, что $\hat x=1$. Соответственно, при $\hat k=1$  получаем, что $\hat x=0$.

Если $x\ne\hat x$, то из соотношений~\eqref{eq:f1}--\eqref{eq:fm} следует, что существуют правые и левые пределы
функции
$f$ --- соответственно $f(x+0)$ и $f(x-0)$. Причем $f(x+0)\ne f(x-0)$ не более чем в сч\"етном числе точек.
Этими точками могут быть только концы отрезков $m$-го порядка при
$m=1,2,\ldots$. Но не все концы указанных отрезков обязательно являются точками разрыва функции $f$.

Укажем тривиальный случай, при котором функция, формально принадлежащая классу $D_1$, принимает конечное число
значений. А именно, это возможно, только если
$\beta_1=\beta_2=\ldots=\beta_{\hat k-1}=\beta_{\hat k+1}=\ldots=\beta_n=d_{\hat k}\beta_1+\beta_{\hat k}$, что следует
из соотношения для $f_2$ в формуле~\eqref{eq:f2}. Самоподобная функция в этой ситуации тождественно равна константе
$\beta_1$.

В частности, при $n=2$ (без ограничения общности считаем, что $\hat k=2$), эти условия сводятся к равенству
$d_{2}\beta_1+\beta_2=\beta_1$.
\end{proof}

\begin{rem}\label{rem:invariant_g}
Всегда будем предполагать,  что при $i,i+1\ne \hat k$ выполнено $\beta_i\ne\beta_{i+1}$. Иначе, можно уменьшить параметр $n$ на единицу,
объединив преобразования $\hat G_i$ и $\hat G_{i+1}$ в одно новое преобразование $\hat G_i$. Из формул~\eqref{eq:f1}--\eqref{eq:fm} следует, что
для при $i\ne \hat k$ операторы $G_i$ и $\tilde G_i$ действуют одинаково.
\end{rem}

В окрестности точки $\hat x$ поведение самоподобной функции $f$ может быть более сложным. Опираясь на
формулы~\eqref{eq:fm} разберём эти случаи подробнее.

Введём для отрезков  $I_{\tiny\underbrace{\hat k,\hat k,\ldots,\hat k}_{m-1}, k_m}$ $m=1,2,\ldots n$ следующее
обозначение $I(m,k_m):=I_{\tiny\underbrace{\hat k,\hat k,\ldots,\hat k}_{m-1},k_m}$.

Из соотношений~\eqref{eq:fm} следует, что для  $i\ne \hat k$ справедливы соотношения
$$
\left. f\right|_{I(m,i)}=d_{\hat k}^{m-1}
\left(\frac{\beta_{\hat k}}{d_{\hat k}-1}+\beta_{i}\right)-\frac{\beta_{\hat k}}{d_{\hat k}-1},\quad \text{ при }
d_{\hat k}\ne 1
$$
и
\begin{gather*}
\left. f\right|_{I(m, i)}=(m-1)\beta_{\hat k}+\beta_{i}\quad\text{ при}\quad d_{\hat k}=1.
\end{gather*}

Из этих формул следует, что в окрестности точки $\hat x$ возможно следующее поведение функции $f$.

\textbf{1.} $|d_{\hat k}|<1$. В этом случае функция $f$ принадлежит пространству $B[0,1]$ (ограниченных на отрезке
функций). Доказательство этого факта полностью аналогично доказательству теоремы~2.1 работы~\cite{Sheip}. При этом,
существует предел значений функции $f$ в точке $\hat x$ и
$$
\lim_{x\to\hat x}f(x)=\dfrac{\beta_{\hat k}}{1-d_{\hat k}}.
$$

Пример такой функции представлен на рис.~1.

\textbf{2.} $d_{\hat k}=1$. Здесь возможны три случая (исключая тривиальный, указанный в доказательстве
теоремы~\ref{tm:schet}).

а) Если $\beta_{\hat k}=0$, $n=3$, $\hat k=2$, то функция $f$ принимает конечное число значений и точка $\hat x$
является точкой разрыва 1-го рода. При этом $f(\hat x-0)=\beta_1$, $f(\hat x+0)=\beta_3$.

б) Если $\beta_{\hat k}=0$ и не ситуация, описанная в 2\;а),  $\hat x$ является точкой разрыва 2-го рода (см. рис.2).

в) Если $\beta_{\hat k}\ne 0$, то
$$
\lim_{x\to\hat x}f(x)=\sign\beta_{\hat k}\cdot(+\infty);
$$

\textbf{3.}  $d_{\hat k}>1$. а) Если все выражения
$\frac{\beta_{\hat k}}{d_{\hat k}-1}+\beta_{j}$, $j\ne \hat k$
имеют один и тот же знак, то существует
$$
\lim_{x\to\hat x} f(x)=\sign \left(\frac{\beta_{\hat k}}{d_{\hat k}-1}+\beta_{j}\right)\cdot(+\infty).
$$

б) Если выражения $\frac{\beta_{\hat k}}{d_{\hat k}-1}+\beta_{j}$ при $j\ne \hat k$ имеют разные знаки, то
$\hat x$ точка разрыва 2-го рода функции $f$.

Нулевое значение в этом случае считается разным и с положительным и с отрицательным значениями.

Заметим, что если $\hat k=2$ и $d_2\beta_1+\beta_2=\beta_1$, то существует $f(\hat x-0)=\beta_1$. Аналогично, если
$\hat k=n-1$ и $d_{n-1}\beta_n+\beta_{n-1}=\beta_n$, то существует $f(\hat x+0)=\beta_n$.

\textbf{4.} Если $d_{\hat k}\leqslant -1$, то $\hat x$ точка разрыва 2-го рода.

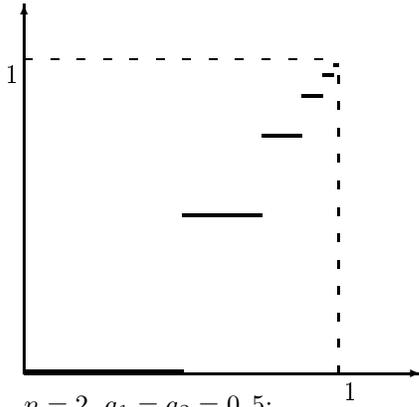
\begin{picture}(360,190)
\put(0,30){
\begin{picture}(180,150)
\put(0,0){\vector(1,0){150}}
\put(0,0){\vector(0,1){140}}

\linethickness{0.4mm}
\put(0,1){\line(1,0){60}}
\put(60,60){\line(1,0){30}}
\put(90,90){\line(1,0){15}}
\put(105,105){\line(1,0){8}}
\put(113,113){\line(1,0){4}}
\put(117,117){\line(1,0){2}}
\thinlines
\multiput(119,0)(0,10){12}{\line(0,1){3}}
\multiput(0,119)(10,0){12}{\line(1,0){3}}
\put(121,-10){\small 1}
\put(-7,110){\small 1}

\vspace{2mm}
\put(0,-15){\small $n=2$, $a_1=a_2=0,5$;}
\put(0,-30){\small $d_1=\beta_1=0$, $d_2=\beta_2=0,5$; $\hat x=1$.}
\put(30,-45){\small Рис. 1 (пример 1)}
\end{picture}}
\put(250,30){
\begin{picture}(180,150)
\put(0,0){\vector(1,0){150}}
\put(0,0){\vector(0,1){140}}

\linethickness{0.4mm}
\put(0,117){\line(1,0){60}}
\put(60,1){\line(1,0){30}}
\put(90,117){\line(1,0){15}}
\put(90,117){\line(1,0){15}}
\put(105,1){\line(1,0){7.5}}
\put(113,117){\line(1,0){4}}

\thinlines
\multiput(119,0)(0,10){12}{\line(0,1){3}}
\put(121,-10){\small 1}
\put(-7,110){\small 1}

\put(0,-18){\small $n=3$, $a_1=0,5$, $a_2=a_3=0,25$;}
\put(0,-33){\small $d_3=\beta_1=1$, $\beta_2=\beta_3=0$; $\hat x=1$.}
\put(30,-48){\small Рис. 2 (пример 2б)}
\end{picture}}
\end{picture}

\begin{picture}(360,200)
\put(0,30){
\begin{picture}(180,150)
\put(0,0){\vector(1,0){150}}
\put(0,0){\vector(0,1){140}}

\linethickness{0.4mm}
\put(0,5){\line(1,0){45}}
\put(90,0){\line(1,0){45}}

\put(45,25){\line(1,0){15}}
\put(75,15){\line(1,0){15}}

\put(60,45){\line(1,0){5}}
\put(70,25){\line(1,0){5}}

\put(65,75){\line(1,0){1.6}}
\put(68.4,45){\line(1,0){1.6}}

\put(66,105){\line(1,0){1}}

\thinlines
\put(135,-1){\line(0,1){3}}
\put(137,-10){\small 1}

\put(0,-15){\small $n=3$, $a_1=a_2=a_3=1/3$;}
\put(0,-30){\small $\beta_1=\beta_2=1$, $d_2=2$, $\beta_3=0$; $\hat x=0,5$.}
\put(30,-45){\small Рис. 3 (пример 3а))}
\end{picture}}
\put(250,30){
\begin{picture}(180,150)
\put(0,60){\vector(1,0){150}}
\put(0,0){\vector(0,1){140}}

\linethickness{0.4mm}
\put(0,65){\line(1,0){45}}
\put(90,50){\line(1,0){45}}

\put(45,75){\line(1,0){15}}
\put(75,35){\line(1,0){15}}

\put(60,85){\line(1,0){5}}
\put(70,15){\line(1,0){5}}

\put(65,105){\line(1,0){1.6}}
\put(68.4,0){\line(1,0){1.6}}

\put(66,105){\line(1,0){1}}

\thinlines
\put(135,59){\line(0,1){3}}
\put(137,50){\small 1}

\put(0,-15){\small $n=3$, $a_1=a_2=a_3=1/3$;}
\put(0,-30){\small $\beta_1=\beta_2=\beta_1=1$, $d_2=2$, $\beta_3=-2$, $\hat x=0,5$.}
\put(30,-45){\small  Рис. 4 (пример 3б))}
\end{picture}}
\end{picture}

\vspace{10mm}

\begin{dif}
Самоподобной струной Стилтьеса будем называть струну, колебания которой описываются задачей~\eqref{eq:vved1}--\eqref{eq:vved2}, где вес $\rho$
представляет собой обобщенную производную самоподобной функции класса $D_1$.
\end{dif}

\begin{rem}\label{rem:indef_string}
Из соотношения~\eqref{eq:osob_tochka} следует, что в отличии от классической стилтьесовской струны, точка накопления масс не обязательно
является концом нити. Кроме того, в этом случае задача может быть индефинитной, т.е. "<массы"> могут принимать отрицательное значение.
\end{rem}

В связи с этим представляет интерес вопрос о нахождении условий на параметры самоподобия, при которых они порождают
неубывающую самоподобную функцию. Для произвольной самоподобной функции, являющейся неподвижной точкой оператора $G$,
определённого формулой~\eqref{eq:oper_pod_ck}, в работе~\cite{Sheip} получены достаточные условия неубывания функции
$f$. Заметим, что в этой работе в формулировке теоремы о достаточных условиях неубывания самоподобной функции допущена
ошибка. Для восстановления правильного утверждения и для полноты изложения приведём формулировку этого результата (см.
теорему 4.1~\cite{Sheip}).

\begin{tm}
Пусть непрерывная слева функция $f$ является  неподвижной точкой сжимающего оператора $G$, определённого формулой~\eqref{eq:oper_pod_ck}, а
параметры самоподобия $\{c_k\}$, $\{d_k\}$ и $\{\beta_k\}$ выбраны так, что $f(0)=1$ и $f(1)=1$. Тогда для того чтобы функция $f$ была
неубывающей необходимо, чтобы при всех $k=1,2,\ldots,n$ выполнялись неравенства
\begin{equation}\label{eq:neobhod_ck}
\beta_{k}\leqslant \beta_{k+1},\quad c_k+d_k+\beta_k\leqslant \beta_{k+1}, \quad c_k+d_k\geqslant 0,
\end{equation}
и достаточно, чтобы при всех $k=1,2,\ldots,n$ выполнялись неравенства
\begin{equation}\label{eq:dostat_ck}
\beta_{k}\leqslant\beta_{k+1}, \quad c_k+d_k+\beta_k\leqslant \beta_{k+1}, \quad c_k\geqslant 0, \quad d_k\geqslant 0,
\end{equation}
\end{tm}
В общей ситуации между необходимыми и достаточными условиями есть "<зазор">. В работе~\cite{Sheip} приведён пример самоподобной функции,
удовлетворяющей условиям~\eqref{eq:neobhod_ck}, но не являющейся монотонной.

Для самоподобной функции нулевого спектрального порядка можно сформулировать критерий монотонности в терминах параметров самоподобия.
\begin{tm}\label{tm:def_string}
Самоподобная функция $f\in D_1$ является неубывающей тогда и только тогда когда её параметры самоподобия удовлетворяют
условиям:

1) $d_{\hat k}>0$;

2) $\beta_1\leqslant\ldots\leqslant \beta_{\hat k-1} \leqslant d_{\hat k}\beta_1+\beta_{\hat k}$
$\leqslant d_{\hat k}\beta_n+\beta_{\hat k}\leqslant\beta_{\hat k+1}\ldots\leqslant\beta_{n}$
при $1<\hat k <n$.

При $\hat k=n$ условие 2) меняется на неравенства:
$\beta_1\leqslant\ldots\leqslant \beta_{n-1}\leqslant \beta_{n-1}\leqslant d_n\beta_1+\beta_n$.

При $\hat k=1$ условие 2) меняется на неравенства:
$d_1\beta_n+\beta_1\leqslant\beta_2\leqslant\ldots\leqslant \beta_{n}$.

Чтобы функция $f$ не возрастала неравенства 2) надо поменять на противоположные.
\end{tm}
\begin{proof}
Рассмотрим случай $1<\hat k<n$. Из формулы\eqref{eq:f1} следует, что для неубывания самоподобной функции на отрезках
первого порядка необходимо и достаточно, чтобы выполнялись неравенства
$\beta_{i}\leqslant\beta_j$ $\forall i<j$, $i,j\ne\hat k$.
Чтобы функция не убывала на отрезках второго порядка необходимо и достаточно, чтобы она не убывала при переходе с
отрезка первого порядка с номером $\hat k -1$ на первый отрезок второго порядка с номером $\hat k, 1$  и при переходе с
последнего отрезка второго порядка с номером $\hat k, n$ на отрезок первого порядка с номером $\hat k +1$. Из формул
\eqref{eq:f2} следует, что для этого необходимо и достаточно выполнения неравенств
$$
\beta_{\hat k-1} \leqslant d_{\hat k}\beta_1+\beta_{\hat k},\qquad
d_{\hat k}\beta_n+\beta_{\hat k}\leqslant\beta_{\hat k+1}.
$$

В силу самоподобия функции $f$ условия неубывания будут выполняться и на отрезках более высокого порядка, что
непосредственно следует из неравенств\eqref{eq:f2}--\eqref{eq:fm}.

Необходимость условие $d_{\hat k}>0$ очевидна.

Случаи $\hat k=1$ и $\hat k=n$ рассматриваются аналогично.
\end{proof}

Легко видеть, что для самоподобной функции $f$ нулевого спектрального порядка, удовлетворяющей условиям $f(0)=0$, $f(1)=1$
условия~\eqref{eq:neobhod_ck}, \eqref{eq:dostat_ck} превращаются в необходимые и достаточные условия монотонности, описанные в
теореме~\ref{tm:def_string}.

\section{Самоподобные функции нулевого спектрального порядка, задаваемые оператором подобия, меняющего ориентацию отрезка}
Как сказано в замечании~\ref{rem:invariant_g} функция, являющаяся неподвижной точкой оператора $G$, не изменится, если операторы $G_k$ поменять
на операторы $\tilde G_k$ при $k\ne\hat k$. А именно, неподвижная функция отображения $G$ принимает на подотрезках $[\alpha_k,\alpha_{k+1}]$
значения, равные $\beta_k$ при $k\ne \hat k$ независимо от того, отображения $\tilde G_k$ или $G_k$ определяют оператор $G$. Поэтому имеет смысл
рассматривать отображение $\tilde G_k$ только при $k=\hat k$.

\subsection{Формулы для приближений самоподобной функции}
Используя свойство оператора $\tilde G_{\hat k}$ менять ориентацию отрезка $[\alpha_{\hat k},\alpha_{\hat k+1}]$, несложно получить формулы для
$m$-ых приближений самоподобной функции, аналогичные формулы~\ref{eq:f1}--\ref{eq:fm}. Начальная функция также принимается равной тождественно
нулю.
\begin{align}
\label{eq:f_orient1}
&f_1(x)=\beta_{k_1}\text{ при } x\in I_{k_1},\quad k_1=1,2,\ldots,n,\\
\label{eq:f_orient2}
&f_2(x)=\left\{\begin{aligned}&\beta_{k_1}, \quad x\in I_{k_1}, \quad k_1\ne \hat k,\\
&d_{\hat k}\beta_{n-k_2+1}+\beta_{\hat k},\quad x\in I_{\hat k,k_2},\quad k_2=1,2,\ldots, n,
\end{aligned}\right.\\
&f_3(x)=\left\{\begin{aligned}&\beta_{k_1}, \quad x\in I_{k_1}, \quad k_1\ne \hat k,\\
&d_{\hat k}\beta_{n-k_2+1}+\beta_{\hat k},\quad x\in I_{\hat k,k_2},\quad k_2\ne\hat k,\\
&d_{\hat k}(d_{\hat k}\beta_{k_3}+\beta_{\hat k})+\beta_{\hat k},\quad x\in I_{\hat k,\hat k,k_3},
\quad k_3=1,2,\ldots,n\\
\end{aligned}\right.\\
\notag
&\ldots\\
\label{eq:f_orientm}
&f_m(x)=\left\{\begin{aligned}& \beta_{k_1}, \quad x\in I_{k_1}, \quad k_1\ne \hat k,\\
&d_{\hat k}\beta_{n-k_2+1}+\beta_{\hat k},\quad x\in I_{\hat k,k_2},\quad k_2\ne \hat k,\\
&\ldots\\
&d_{\hat k}^{m-2}\beta_{k_{m-1}}+\beta_{\hat k}(1+d_{\hat k}+\ldots+d_{\hat k}^{m-3}),\quad
x\in I_{\tiny\underbrace{\hat k,\hat k,\ldots,\hat k}_{m-2},k_{m-1}},\quad k_{m-1}\ne \hat k\\
&d_{\hat k}^{m-1}\beta_{n-k_m+1}+\beta_{\hat k}(1+d_{\hat k}+\ldots+d_{\hat k}^{m-2}),\quad x\in I_{\tiny\underbrace{\hat k,\hat k,\ldots,\hat
k}_{m-1},k_m},
\quad k_m=1,2,\ldots, n,
\end{aligned}\right.
\end{align}
при чётном $m$. При нечётном $m$ формула~\ref{eq:f_orientm} принимает вид
$$
f_m(x)=\left\{\begin{aligned}& \beta_{k_1}, \quad x\in I_{k_1}, \quad k_1\ne \hat k,\\
&d_{\hat k}\beta_{n-k_2+1}+\beta_{\hat k},\quad x\in I_{\hat k,k_2},\quad k_2\ne \hat k,\\
&\ldots\\
&d_{\hat k}^{m-2}\beta_{n-k_{m-1}+1}+\beta_{\hat k}(1+d_{\hat k}+\ldots+d_{\hat k}^{m-3}),\quad
x\in I_{\tiny\underbrace{\hat k,\hat k,\ldots,\hat k}_{m-2},k_{m-1}},\quad k_{m-1}\ne \hat k\\
&d_{\hat k}^{m-1}\beta_{k_m}+\beta_{\hat k}(1+d_{\hat k}+\ldots+d_{\hat k}^{m-2}),\quad x\in I_{\tiny\underbrace{\hat k,\hat k,\ldots,\hat
k}_{m-1},k_m},
\quad k_m=1,2,\ldots, n.
\end{aligned}\right.
$$

\subsection{Координаты особой точки}
Рассмотрим, как меняются формулы для определения координат особой точки, если преобразование c номером $\hat k$ меняет ориентацию отрезка
$[\alpha_{\hat k},\alpha_{\hat k+1}]$. Координаты левого и правого концов отрезка $I_{\underbrace{\hat k,\hat k,\ldots,\hat k}_{m}}$ в этом случае имеют вид
\begin{gather*}
\tilde x_{l}=\alpha_{\hat k}+a_{\hat k}(1-\alpha_{\hat k+1})+a^2_{\hat k}\alpha_{\hat k}+a_{\hat k}^2(1-\alpha_{\hat k+1})+\ldots+a^{\lceil m/2\rceil}_{\hat k}A,\\
\tilde x_{r}=1-(1-\alpha_{\hat k+1})-a_{\hat k}\alpha_{\hat k}-a^2_{\hat k}(1-\alpha_{\hat k+1})-a_{\hat k}^2\alpha_{\hat k}-\ldots-
a^{\lceil m/2\rceil}_{\hat k}A,
\end{gather*}
где
$$
A=\begin{cases}\alpha_{\hat k} & \text{при } m=2l-1\\
1-\alpha_{\hat k+1} & \text{при } m=2l\end{cases},
$$
а $\lceil x\rceil$ --- наименьшее целое число, не меньшее $x$ ("<потолок">), $l=1,2,\ldots$.

Несложно вычислить, что особая точка $\hat x$ в этом случае имеет координату
\begin{equation}\label{eq:osob_tochka2}
\hat x=\dfrac{\alpha_{\hat k+1}}{1+a_{\hat k}}.
\end{equation}

Отметим, что в этом случае, независимо от значения $\hat k$ особая точка $\hat x$ не может совпасть с концами отрезка $[0,1]$.

\subsection{Условия неубывания функции в случае $\tilde G_{\hat k}$.}
На основании формул~\ref{eq:f_orient1}--\ref{eq:f_orientm} несложно получить аналог теоремы~\ref{tm:def_string}.

\begin{tm}
Самоподобная функция $f\in D_1$, заданная таким оператором подобия,  у которого $\tilde G_{\hat k}$ меняет ориентацию отрезка
$[\alpha_{\hat k},\alpha_{\hat k+1}]$, является неубывающей тогда и только тогда когда её параметры самоподобия удовлетворяют условиям:

1) $d_{\hat k}<0$;

2) $\beta_1\leqslant\beta_2\leqslant\ldots\leqslant \beta_{\hat k-1} \leqslant d_{\hat k}\beta_n+\beta_{\hat k}$
$\leqslant d_{\hat k}\beta_1+\beta_{\hat k}\leqslant\beta_{\hat k+1}\ldots\leqslant\beta_{n}$
при $1<\hat k <n$.

При $\hat k=n$ условие 2) примет вид: $\beta_1\leqslant\beta_2\leqslant\ldots\leqslant \beta_{n-1}\leqslant \beta_{n}$ и
$\beta_{n-1} \leqslant d_{n}\beta_n+\beta_{n}$.

При $\hat k=1$ условие 2) примет вид: $\beta_1\leqslant\beta_2\leqslant\ldots\leqslant \beta_{n-1}\leqslant \beta_{n}$ и
$d_1\beta_1+\beta_1\leqslant\beta_2$.
\end{tm}

\begin{rem}
Функцию нулевого спектрального порядка, у которой преобразование $\tilde G_{\hat k}$ меняет ориентацию отрезка
$[\alpha_{\hat k},\alpha_{\hat k+1}]$ можно задать другим набором параметров самоподобия, для которых все преобразования $G_i$
сохраняют ориентацию отрезков.
\end{rem}

Чтобы найти эти параметры самоподобия, достаточно заметить, что сжимающее отображение $G^2$ имеет ту же неподвижную функцию $f$, что и
преобразование $G$.  Введём преобразование $F:=G^2$. Пусть $N$,
$\alpha(F)_k$,
$\beta(F)_k$ и
$d(F)_k$,
$k=1,2,\ldots, N$
--- его параметры самоподобия. Через $m$ обозначим тот единственный индекс, для которого $d(F)_m\ne 0$. Эти числа можно вычислить через
параметры самоподобия исходного преобразования $G$ следующим образом:
\begin{gather*}
 N=2n-1, \quad m=n,\quad  d(F)_{m}=d_{\hat k}^2,\\
 a(F)_k=a_k, \quad k=1,2,\ldots, \hat k-1,\\
 a(F)_{m+k}=a_{\hat k}\cdot a_{n-k}, \quad k=0,1,\ldots,n-1,\\
 a(F)_{m+n+k}=a_{\hat k+k+1}, \quad k=0,1,\ldots,n-\hat k,\\
 \beta(F)_k=\beta_k, \quad k=1,2,\ldots, \hat k-1,\\
 \beta(F)_{m+k}=d_{\hat k}\cdot \beta_{n-k}, \quad k=0,1,\ldots,n-1,\\
 \beta_{m+n+k}=\beta_{\hat k+k+1}, \quad k=0,1,\ldots,n-\hat k,
\end{gather*}

Как и преобразование $G$, gреобразование $F$ также можно задать с помощью отображений $F_k$, сопоставляющих функции
$f\in L_p[0,1]$ функции $F_k(f)\in L_p[\alpha(F)_k,\alpha(F)_{k+1}]$, $k=1,2,\ldots, N$:
$$
(F_k f)(x)=f(t), \quad x=a(F)_kt+\alpha(F)_k, \quad x\in[\alpha(F)_k,\alpha(F)_{k+1}], \quad t\in[0,1].
$$

 У преобразования $F$ все его составляющие преобразования $F_k$, $k=1,2,\ldots, N$ сохраняют ориентацию отрезков.

Автор признателен А.И.~Назарову за полезные обсуждения и замечания.


\begin{thebibliography}{99}


\bibitem{KacKr} И.~С.~Кац, М.~Г.~Крейн. \textit{О спектральных функциях
струны.} В кн.:~Ф.~Аткинсон. \textit{Дискретные и непрерывные граничные задачи}, М., ", 1968,
стр.~648--733.

\bibitem{SV} M.~Solomyak, E.~Verbitsky. \textit{On a spectral problem
related to self-similar measures}//Bull. London Math. Soc.,
\textbf{27}\,(1995), pp.~242--248.

\bibitem{VlSh1} А.~А.~Владимиров, И.~А.~Шейпак, \emph{Самоподобные функции в пространстве \(L_2[0,1]\) и задача
Штурма-Лиувилля с сингулярным индефинитным весом}//Математический сборник, \textbf{197}(11), 2006, с. 13-30

\bibitem{VlSh2} А.~А.~Владимиров, И.~А.~Шейпак, \emph{Индефинитная задача Штурма-Лиувилля для некоторых классов
самоподобных сингулярных весов}/Труды МИРАН, т.255, 2006, с. 88-98

\bibitem{VlSh3} А.~А.~Владимиров, И.~А.~Шейпак, \emph{Асимптотика собственных значений задачи
Штурма-Лиувилля с дискретным самоподобным весом}//http://www.arxiv.org/math/arXiv:0709.0424v1

\bibitem{Bar} M.~Barnsley, \emph{Fractals everywere}//Academic Press, 1988.

\bibitem{deRam1} De Rham G., \emph{Une peu de math\'ematique \`a propos d'une courbe plane}// Rev. de
math. elemetaires, \textbf{2}:4,5 (1947), pp.73--76, 89--97.

\bibitem{deRam2} De Rham G., \emph{Sur une courbe plane}// J. Math. Pures Appl. (9), \textbf{35}
(1956), pp. 25--42.

\bibitem{deRam3} De Rham G., \emph{Sur les courbes limit\'es de polygones obtenus par
trisection}//Enseign. Math., (2), \textbf{5} (1959), pp.29--43.

\bibitem{Prot1} Protasov V., \emph{Refinement equations with nonnegative coefficients}//J. Fourier
Anal. Appl., \textbf{6}:1, (2000), pp. 55--78.

\bibitem{Prot2} В.~Ю.~Протасов, \emph{Фрактальные кривые и всплески}//Изв.РАН. Серия матем.,
\textbf{70}:5 (2006), стр.105--145.

\bibitem{Sheip} И.~А.~Шейпак, \emph{О конструкции и некоторых свойствах самоподобных функций в пространствах
$L_p[0,1]$}//Матем. заметки, \textbf{81}:6, (2007), с. 924-938.

\bibitem{ANaz}А.~И.~Назаров, {\em Логарифмическая асимптотика малых уклонений для некоторых гауссовских процессов в $L_2$-норме относительно самоподобной
меры}, Вероятность и статистика. ЗНС ПОМИ, {\bf 311} (2004), 190--213.


\end{thebibliography}
\end{document}